\documentclass[final,leqno,onefignum,onetabnum]{siamltex1213}
\usepackage{amsmath}
\usepackage{amssymb}
\usepackage{amsfonts}
\newcommand{\ppp}{\partial}
\newcommand{\ooo}{\overline}
\def\va{\varphi}
\def\vadd{\varphi_{d_1}}
\def\vad{\varphi_d}
\def\psid{\psi_d}

\def\la{\lambda}
\def\ppp{\partial}
\def\ooo{\overline}
\def\OOO{\Omega}
\def\ddd{\mbox{div}\thinspace}

\def\weight{e^{2s\va}}
\def\o{\Omega}

\def\la{\lambda}
\def\R{\Bbb R}

\title{Stability for some inverse problems for transport equations}
\author{Fikret G\"olgeleyen\footnotemark[2]
\and Masahiro Yamamoto\footnotemark[3]}

\begin{document}
\maketitle
\slugger{mms}{xxxx}{xx}{x}{x--x}
\renewcommand{\thefootnote}{\fnsymbol{footnote}}
\footnotetext[2]{Department of Mathematics, Bulent Ecevit University, Zonguldak, 67100 Turkey (f.golgeleyen@beun.edu.tr). The work has been done during the stay of the first author
at Graduate School of Mathematical Sciences of The University of Tokyo, which was supported by Leading Graduate Course for Frontiers of Mathematical Sciences and Physics.}
\footnotetext[3]{Department of Mathematical Sciences, The University of Tokyo, 3-8-1 Komaba, Meguro, Tokyo, 153-8914 Japan (myama@ms.u-tokyo.ac.jp). This author is partially supported by Grant-in-Aid for Scientific
Research (S) 15H05740 of Japan Society for the Promotion of Science.}

\begin{abstract}
In this article, we consider inverse problems of determining a source term and a coefficient of
a first-order partial differential equation and prove conditional stability estimates with minimum boundary observation data
and relaxed condition on the principal part.
\end{abstract}



\pagestyle{myheadings}
\thispagestyle{plain}
\markboth{FIKRET G\"OLGELEYEN AND MASAHIRO YAMAMOTO}{STABILITY FOR SOME INVERSE PROBLEMS FOR TRANSPORT EQUATIONS}

\section{Introduction and main results}

Let $\Omega\subset {\Bbb R}^n$ be a bounded domain with smooth
boundary $\ppp\Omega$ and let $\nu(x)$ be the unit outward normal vector to
$\ppp\Omega$.  Let us consider
$$
\ppp_ty(x,t) + H(x)\cdot \nabla y(x,t) + V(x)y(x,t) = f(x)R(x,t),
\quad  x \in \Omega, \thinspace 0<t<T         \eqno{(1.1)}
$$
and
$$
y(x,0) = 0, \qquad x \in \Omega.        \eqno{(1.2)}
$$
We assume that $H := (h_1, ..., h_n) \in \{C^1(\ooo\o)\}^n$ and
$V \in L^{\infty}(\Omega)$.

Throughout this paper, we set $x = (x_1, ..., x_n) \in {\Bbb R}^n$,
$\ppp_j = \frac{\ppp}{\ppp x_j}$ for $j=1,2,..., n$ and
$\ppp_t = \frac{\ppp}{\ppp t}$, $\nabla = (\ppp_1, ..., \ppp_n)$,
$\Delta = \sum_{j=1}^n \ppp_j^2$, and $H\cdot J$ denotes the scalar product
of $H, J \in {\R}^n$.
\\

The main problems in this paper are
\\
\vspace{0.2cm}
{\bf Inverse source problem}\\
Let $H$, $V$, $R$, $\Gamma \subset \ppp\Omega$, $T>0$ be given suitably.
Determine $f(x)$, $x \in \Omega$ from $y\vert_{\Gamma\times (0,T)}$.
\\

Moreover we consider
$$
\ppp_tu(x,t) + H(x)\cdot \nabla u(x,t) + V(x)u(x,t) = 0,
\quad  x \in \Omega, \thinspace 0<t<T         \eqno{(1.3)}
$$
and
$$
u(x,0) = a(x), \qquad x \in \Omega.        \eqno{(1.4)}
$$
{\bf Inverse coefficient problem}\\
Let $a$ and $H$ be suitably given.  Determine $V(x)$ and/or $H(x)$ by
data $u\vert_{\Gamma\times (0,T)}$.

Equations (1.1) and (1.3) are transport equations and are models
in physical phenomena such as Liouville equation and
the mass conservation law.  Moreover the transport equation is related to
the integral geometry (e.g., Amirov [1]).
As for other physical backgrounds such as neutron transport and medical
tomography, see e.g., Case and Zweifel [8], Ren, Bal and Hielscher [18].

Our inverse problem is formulated with a single measurement, and
Gaitan and Ouzzane [9], Klibanov and Pamyatnykh [15],
Machida and Yamamoto [17] discuss the uniqueness and the stability for
inverse problems for initial/boundary value problems for transport equations
by Carleman estimates.
The papers [15] and [17] discuss transport
equations with integral terms where solutions $y$ and $u$ depend also on the
velocity as well as the location $x$ and the time $t$.

The main methodology in [9], [15], [17] is based on Bukhgeim and Klibanov [7].
After that, there have been many works.  Limited to hyperbolic and parabolic
equations, we can refer for example to Baudouin, de Buhan and Ervedoza
[3], Beilina and Klibanov [4], Bellassoued and Yamamoto [6],
 Imanuvilov and Yamamoto [11], [12],
Klibanov [14], Yamamoto [21] and the references therein.
Here we do not intend to give any complete lists of the references.
In [9] and [15], the key Carleman estimate is the same as
the Carleman estimate for a second-order hyperbolic equation
and in order to apply the Carleman estimate one has to extend the
solutions $y$ and $u$ to (1.1) and (1.3) to the time interval
$(-T,0)$.  Such an extension argument makes the proofs longer, and requires
an extra condition to unknown coefficients and initial value as in [15].
In Sections 2 and 4, we prove Carleman estimates (Lemmata 1 and 3),
which can directly estimate initial values.
Thanks to our Carleman estimates, we can
simplify the proofs of the stability and relax the constraints
of the principal coefficients $H$'s.

As for inverse problems for transport equations with many measurements,
see surveys Bal [2], Stefanov [20]
and the references therein.  Klibanov and Yamamoto [16] established the
exact controllability for the transport equations by a Carleman estimate.

Unlike [2] and [20], we discuss the inverse problems for a single
initial/boundary value problem where we need not change initial values or
boundary values.
More precisely, in the formulation for the inverse problems in [2] and [20],
we have to change boundary inputs on some subboundary and repeat
measurements of the corresponding boundary outputs on  other subboundary.
One can apply the method of characteristics to the same kind of inverse problem
for the first-order equation and see Belinskij [5],
Chapter 5 of Romanov [19] for example.

\vspace{0.2cm}

We set
$$
Q = \OOO\times (0,T)
$$
and
$$
\left\{ \begin{array}{rl}
&\ppp\Omega_+ = \{ x\in \ppp\Omega; \thinspace (\nu(x)\cdot H(x)) > 0\},\\
&\ppp\Omega_- = \{ x\in \ppp\Omega; \thinspace (\nu(x)\cdot H(x)) < 0\}.\\
\end{array}\right.
$$

Throughout this paper, we assume that $\psi \in C^2(\ooo\Omega)$
and $H = (h_1, ..., h_n) \in \{C^1(\ooo\Omega)\}^n$ satisfy
$$
\mu:= \min_{x\in\ooo\OOO} (H(x)\cdot \nabla \psi(x)) > 0.
                                            \eqno{(1.5)}
$$
We here note by (1.5) that $\vert H(x)\vert \ne 0$ for
$x \in \ooo\Omega$.

We give four cases where (1.5) holds.
\\
{\bf Case 1}.  We assume
$$
\vert \nabla d\vert > 0 \quad \mbox{on $\ooo\OOO$}, \qquad
H(x) = \nabla d(x), \qquad x \in \ooo\OOO
$$
with some $d\in C^2(\ooo\OOO)$.
Then (1.5) holds if we choose $\psi(x) = d(x)$,
$x \in \OOO$.
\\
{\bf Case 2.}
Let us assume that $\{ (h_1(x), ..., h_n(x));
\thinspace x\in \ooo\o\} \subset {\Bbb R}^n$ is
separated from $(0,..., 0)$ by a hyperplane
$a_1x_1 + \cdots + a_nx_n = 0$ with some $a_1, ..., a_n \in {\Bbb R}$
and $\vert a_1\vert + \cdots + \vert a_n\vert \ne 0$.
Then $\psi(x) = a_1x_1 + \cdots + a_nx_n$ or
$\psi(x) = -a_1x_1 - \cdots - a_nx_n$ satisfies (1.5).  In particular,
(1.5) holds if $H(x)$ is a constant vector because
$\{ H(x); \thinspace \ooo\o\}$ is composed of one point.
In fact, the separation condition means that
$\vert (H(x)\cdot\nabla\psi(x)) \vert
= \vert a_1h_1(x) + \cdots
+ a_nh_n(x) \vert > 0$ for all $x\in \ooo\o$ or $< 0$ for all $x \in
\ooo\o$.
\\
{\bf Case 3.}
Let $0 \in \o$.  We assume that
there exists a constant $\delta_0 > 0$ such that
$$
\vert H(x)\vert \ge \delta_0, \quad x \in \ooo\OOO.
$$
Then $\psi(x) = \sum_{j=1}^n x_jh_j(x)$ satisfies (1.5) if
$\max_{x\in\ooo\o} \vert x\vert$ is sufficiently small.
\\
{\bf Proof.} By the Cauchy-Schwarz inequality, we have
\begin{eqnarray*}
&& (H(x)\cdot \nabla\psi(x))
= \sum_{\ell=1}^n h_{\ell}(x)^2 + \sum_{\ell=1}^nh_{\ell}(x)
\sum_{j=1}^n x_j\ppp_{\ell}h_j(x)\\
\ge &&\min_{x\in\ooo\o} \vert H(x)\vert^2
- \left(\sum_{\ell=1}^n h_{\ell}(x)^2\right)^{\frac{1}{2}}
\left( \sum_{\ell=1}^n \left\vert
\sum_{j=1}^n x_j\ppp_{\ell}h_j(x)\right\vert^2\right)^{\frac{1}{2}}\\
\ge && \delta_0^2 - \Vert H\Vert_{\{L^{\infty}(\o)\}^n}
\left(\sum_{\ell=1}^n \left( \sum_{j=1}^n \vert x_j\vert^2\right)
\left( \sum_{j=1}^n \vert \ppp_{\ell}h_j(x)\vert^2\right)
\right)^{\frac{1}{2}}\\
\ge&& \delta_0^2 -  \max_{x\in\ooo\o} \vert x\vert
\Vert H\Vert_{\{L^{\infty}(\o)\}^n}
\Vert \nabla H\Vert_{\{L^{\infty}(\o)\}^{n\times n}}.
\end{eqnarray*}
Therefore (1.5) holds true if
$$
0 < \max_{x\in\ooo\OOO} \vert x\vert
< \frac{\min_{x\in \ooo\OOO} \vert H(x)\vert^2}
{\Vert H\Vert_{\{L^{\infty}(\OOO)\}^n}
\Vert \nabla H\Vert_{\{L^{\infty}(\o)\}^{n\times n}}}.
$$
{\bf Case 4}.  We assume that there exists $i_0 \in \{1, 2, ..., n\}$ such that
$h_{i_0}(x) > 0$ for all $x \in \ooo\OOO$.
Then we choose sufficienly small $b \in \R$ such that
$\ooo\OOO \subset \{ (x_1, x_2, ..., x_n);\thinspace
x_{i_0} > b\}$.  Setting $\psi(x) = (x_{i_0}-b)^2$, we verify that
(1.5) holds.  In fact,
$$
(H(x)\cdot \nabla \psi(x)) = 2h_{i_0}(x)(x_{i_0} - b) > 0 \quad
\mbox{for $x \in \ooo\OOO$}.
$$
\vspace{0.2cm}

Now we state the first main result concerning the stability for the inverse
source problem.
\\
{\bf Theorem 1}\\
Let $y \in H^1(Q)$ satisfy (1.1) and (1.2), and let (1.5) be
satisfied with some constant $\mu>0$.  We assume that
$$
R(x,0) \ne 0, \quad x \in\ooo\Omega          \eqno{(1.6)}
$$
and
$$
\ppp_ty, \ppp_tR \in H^1(Q), \quad
\ppp_tR \in L^2(0,T;L^{\infty}(\o)).
$$
Let
$$
T > \frac{\max_{x\in\ooo\OOO} \psi(x) - \min_{x\in \ooo\OOO} \psi(x)}
{\mu}.                             \eqno{(1.7)}
$$
\\
(i) We assume
$$
\Vert \ppp_ty\Vert_{L^2(Q)} \le M
$$
with fixed constant $M>0$.  Then there exist constants $\theta \in (0,1)$
and $C>0$, which are dependent on $\Omega, T, H, \Vert V\Vert_{L^{\infty}
(\OOO)}, \psi,
M, \Vert \ppp_tR\Vert_{L^2(0,T;L^{\infty}(\OOO))}$, such that
$$
\Vert f\Vert_{L^2(\OOO)} \le C\left\{\left(\int^T_0\int_{\ppp\Omega_+}
(H\cdot\nu)\vert \ppp_ty\vert^2 dS_xdt\right)^{\frac{\theta}{2}}
+ \left(\int^T_0\int_{\ppp\Omega_+}
(H\cdot\nu)\vert \ppp_ty\vert^2 dS_xdt\right)^{\frac{1}{2}}\right\}
$$
for all $f \in L^2(\o)$.
\\
(ii) Without the assumption in (i) concerning $\Vert\ppp_ty\Vert_{L^2(Q)}$,
there exists a constant $C>0$, which is dependent on $\Omega, T, H,
\Vert V\Vert_{L^{\infty}(\OOO)}, \psi,
\Vert \ppp_tR\Vert_{L^2(0,T;L^{\infty}(\OOO))}$, such that
$$
\Vert f\Vert_{L^2(\OOO)} \le C\left( \int^T_0\int_{\ppp\Omega}
\vert (H\cdot\nu)\vert \vert \ppp_ty\vert^2 dS_xdt\right)^{\frac{1}{2}}
$$
for all $f \in L^2(\OOO)$.
\\
(iii) In addition to (1.1) and (1.2), we assume
$$
y = 0 \qquad \mbox{on $\ppp\OOO_- \times (0,T)$}.   \eqno{(1.8)}
$$
Then there exists a constant $C>0$, which is dependent on $\Omega, T, H,
\Vert V\Vert_{L^{\infty}(\OOO)}$, $\psi$,
$\Vert \ppp_tR\Vert_{L^2(0,T;L^{\infty}(\OOO))}$, such that
$$
C^{-1}\left( \int^T_0\int_{\ppp\Omega_+}
(H\cdot\nu) \vert \ppp_ty\vert^2 dS_xdt\right)^{\frac{1}{2}}
\le \Vert f\Vert_{L^2(\OOO)}
\le C\left( \int^T_0\int_{\ppp\Omega_+}
(H\cdot\nu) \vert \ppp_ty\vert^2 dS_xdt\right)^{\frac{1}{2}}  \eqno{(1.9)}
$$
for all $f \in L^2(\OOO)$.
\\
\vspace{0.2cm}

The conclusion of (i) is a stability estimate of H\"older type
and holds under a priori
boundedness $\Vert \ppp_ty\Vert_{L^2(Q)} \le M$, which is called conditional
stability.  On the other hand, the conclusions of (ii) and (iii) are
Lipschitz stability and in particular, with (1.8) we can have
both-sided estimate (1.9) for our inverse problem.

We apply Theorem 1 to the inverse coefficient problem of determining
$V(x)$.
\\
\vspace{0.2cm}
{\bf Theorem 2}\\
For $j=1,2$, let $u_j, \ppp_tu_j \in H^1(Q)$ and let
$$
\ppp_tu_j + H(x)\cdot \nabla u_j + V_j(x)u_j(x,t) = 0 \quad
\mbox{in $Q$},                       \eqno{(1.10)}
$$
$$
u_j(x,0) = a(x), \qquad x \in \OOO    \eqno{(1.11)}
$$
and
$$
u_j = h(x,t) \qquad \mbox{on $\ppp\OOO_-\times (0,T)$}  \eqno{(1.12)}
$$
with suitably given $a$ and $h$.  We assume that there exists
$\psi \in C^2(\ooo\OOO)$ satisfying (1.5) for $H$, (1.7) holds and
$$
u_j, \ppp_tu_j\in H^1(Q)\cap L^2(0,T;L^{\infty}(\OOO)), \quad
j=1,2.
$$
Moreover we assume
$$
\Vert V_j\Vert_{L^{\infty}(\OOO)}, \Vert \ppp_tu_j\Vert_{L^2(0,T;
L^{\infty}(\OOO))} \le M, \quad j=1,2                       \eqno{(1.13)}
$$
and
$$
\vert a\vert > 0 \qquad \mbox{on $\ooo\OOO$},    \eqno{(1.14)}
$$
where $M>0$ is arbitrarily fixed constant.

Then there exists a constant $C>0$ depending on $\Omega, T, H, \psi, a, M$
such that
$$
C^{-1}\Vert \ppp_t(u_1-u_2)\Vert_{L^2(\ppp\OOO_+\times (0,T))}
\le \Vert V_1-V_2\Vert_{L^2(\OOO)}
\le C\Vert \ppp_t(u_1-u_2)\Vert_{L^2(\ppp\OOO_+\times (0,T))}.
                                                    \eqno{(1.15)}
$$
\\

Similarly to Theorem 2, we can discuss the determination of $H$, but
we have to determine $n$ functions as the components of $H$, and so
repeats of measurements of boundary data after changing initial values
suitably are necessary but arguments can be repeated similarly to
the proof of Theorem 2.  Here we omit detailed discussions for the
determination of all the components of $H$, but we consider the
determination of the potential in the case of potential flows.
That is, we consider
$$
\left\{ \begin{array}{rl}
& \ppp_t\rho(x,t) + \ddd (\rho \nabla d(x)) = 0 \quad
\mbox{in $Q$}, \\
& \rho(x,0) = a(x), \qquad x \in \OOO.
\end{array}\right.
\eqno{(1.16)}
$$
The first-order partial differential equation in (1.16) describes the
mass conservation under
stationary potential flow $\nabla d(x)$.  Then we discuss the inverse
problem of determining the stationary potential $d$.

For the statement of the main result, we define an admissible set of
unknown potentials $d$'s.  Let constants $M>0$, $\delta_0>0$ and functions
$g_1 \in C^2(\ppp\OOO)$, $g_2 \in C^1(\ppp\OOO)$ be arbitrarily
chosen.  We define the admissible set of $d$'s by
$$
\mathcal{D} = \mathcal{D}(\delta_0,M,g_1,g_2) :=
\{ d \in C^2(\ooo\OOO);\thinspace \vert \nabla d\vert \ge \delta_0 > 0
\thinspace \mbox{on $\ooo\OOO$}, \\\
$$
$$
\Vert d\Vert_{C^2(\ooo\OOO)} \le M, \quad
d\vert_{\ppp\OOO} = g_1, \thinspace \ppp_{\nu}d\vert_{\ppp\OOO}
= g_2\}.                          \eqno{(1.17)}
$$
For the function $g_2$ given in (1.17), we set
$$
\left\{\begin{array}{rl}
& \Gamma_+ = \{ x\in \ppp\OOO; \thinspace g_2(x) > 0\},\\
& \Gamma_- = \{ x\in \ppp\OOO; \thinspace g_2(x) < 0\}.\\
\end{array}\right.
\eqno{(1.18)}
$$
We state our final main result.
\\
{\bf Theorem 3}\\
For $j=1,2$, let $\rho_j, \ppp_t\rho_j \in H^1(Q)$ and let
$$
\ppp_t\rho_j(x,t) + \ddd (\rho_j \nabla d_j(x)) = 0 \quad
\mbox{in $Q$}                                \eqno{(1.19)}
$$
and
$$
\rho_j = h \quad \mbox{on $\Gamma_-\times (0,T)$}, \qquad
\rho_j(x,0) = a(x), \quad x \in \OOO           \eqno{(1.20)}
$$
with suitable $a$ and $h$.
We assume (1.14),
$$
T > \frac{\sup_{d\in \mathcal{D}}(\max_{x\in\ooo\OOO} d(x)
- \min_{x\in\ooo\OOO} d(x))}{\delta_0^2}               \eqno{(1.21)}
$$
and
$$\left\{
\begin{array}{rl}
& \ppp_t\rho_j \in L^2(0,T;W^{1,\infty}(\OOO)) \cap H^1(Q), \\
& \Vert \rho_j\Vert_{L^2(0,T;H^1(\OOO))} \le M, \quad j=1,2.
\end{array}\right.
                                  \eqno{(1.22)}
$$
Then there exists a constant $C>0$ depending on $\Omega, T, a,
\delta_0, M, g_1, g_2$ such that
$$
C^{-1}\Vert \ppp_t(\rho_1-\rho_2)\Vert_{L^2(\Gamma_+\times (0,T))}
\le \Vert d_1-d_2\Vert_{H^2(\OOO)}
\le C\Vert \ppp_t(\rho_1-\rho_2)\Vert_{L^2(\Gamma_+\times (0,T))}
                                                    \eqno{(1.23)}
$$
for all $d_1, d_2 \in \mathcal{D}$.
\\
\vspace{0.2cm}

The proofs of Theorems 1-3 are based on an argument by the Carleman estimates,
which was originated by  Bukhgeim and Klibanov [7].
Here we used a modified argument by Imanuvilov and Yamamoto [11], [12]
which discussed for inverse problems for second-order hyperbolic equations.

The paper is composed of four sections and an appendix.
In Section 2, we prove
a relevant Carleman estimate for the proofs of Theorems 1 and 2, and an
energy estimate, and in Section 3, the proofs of Theorems 1 and 2 are
completed.
In Section 4, we prove another Carleman estimate for the proof of
Theorem 3 whose weight function gives a Carleman estimate also for the
Laplacian and complete the proof of Theorem 3.
In Appendix, we prove the Carleman estimate for the Laplacian.
\section{Key Carleman estimate and energy estimate}
We recall that
$$
Q = \OOO\times (0,T)
$$
and we set
$$
Pu = \ppp_tu + H(x)\cdot\nabla u + V(x)u, \quad
P_0u = \ppp_tu + H(x)\cdot\nabla u, \quad (x,t) \in Q,
$$
$$
M_0 = \beta \Vert \mbox{div}\thinspace H\Vert_{L^{\infty}(\OOO)}
\Vert \mbox{div}\thinspace (H(H\cdot \nabla\psi))\Vert
_{L^{\infty}(\Omega)}
$$
and
$$
\va(x,t) = -\beta t + \psi(x), \quad (x,t) \in Q   \eqno{(2.1)}
$$
with $\psi \in C^2(\ooo\o)$ and $\beta > 0$, and
$$
B(x) := \ppp_t\va + (H\cdot \nabla\va) = -\beta + (H(x)\cdot\nabla\psi),
\quad x\in \OOO.                   \eqno{(2.2)}
$$

First we prove
\\
{\bf Lemma 1}\\
(i) We have
\begin{eqnarray*}
&& s\int_{\o} B(x)\vert u(x,0)\vert^2 e^{2s\va(x,0)} dx
+ s^2\int_Q B^2(x)\vert u(x,t)\vert^2e^{2s\va} dxdt\\
\le&& 2\int_Q \vert Pu\vert^2 \weight dxdt
+ (sM_0+2\Vert V\Vert_{L^{\infty}(\OOO)}^2)
\int_Q \vert u(x,t)\vert^2e^{2s\va} dxdt\\
+ && s\int^T_0\int_{\ppp\o} B(x)(\nu\cdot H) \vert u\vert^2 \weight dS_xdt
\end{eqnarray*}
for all $s >0$ and $u \in H^1(Q)$ satisfying $u(\cdot,T) = 0$ in
$\OOO$.
\\
(ii) We assume (1.5) and
$$
0 < \beta < \mu:= \min_{x\in \ooo\o} (H(x) \cdot \nabla\psi(x)). \eqno{(2.3)}
$$
Then
$$
s\int_{\o} \vert u(x,0)\vert^2 e^{2s\va(x,0)} dx
+ \frac{s^2(\mu-\beta)^2}{2}\int_Q \vert u(x,t)\vert^2e^{2s\va} dxdt
$$
$$
\le 2\int_Q \vert Pu\vert^2 \weight dxdt
+ s\int^T_0\int_{\ppp\o_+} B(H\cdot\nu)\vert u\vert^2 \weight dS_xdt
                                  \eqno{(2.4)}
$$
for all $s \ge s_0$ and $u \in H^1(Q)$ satisfying $u(\cdot,T) = 0$ in
$\OOO$, where
$$
s_0 = \max\left\{ \frac{4M_0}{(\mu-\beta)^2},
\frac{\sqrt{8}\Vert V\Vert_{L^{\infty}(\OOO)}}{\mu-\beta} \right\}.
$$
\\

Inequality (2.4) is an estimate of Carleman's type, which holds
uniformly for sufficiently large $s>0$.
We emphasize that Carleman estimate (2.4) can estimate also the initial
value $u(x,0)$, and the weight function is linear in $t$.
In [17] such a linear weight function is used for proving a
stability estimate for an inverse problem
for a transport equation.
For the transport equation, the works [9] and [15] used the same
weight function as the second-order hyperbolic equation, that is,
$\va(x,t) = e^{\lambda(d(x) - \beta t^2)}$.
Our choice (2.1) of the weight function enables us not to need to
take any extensions of $u$ to $(0,-T)$ in discussing the inverse problem.
In [15], the extension argument requires an extra assumption on the
initial value and unknown coefficients in addition to (1.14).\\

\vspace{0.3cm}
{\bf Proof of Lemma 1 (i).}\\
First we assume that $V=0$.
We set $w(x,t)=e^{s\va(x,t)}u(x,t)$ and $(L w)(x,t) = e^{s\va(x,t)}P_0
(e^{-s\va}w)$.
Then
$$
Lw = \{\ppp_tw + (H(x)\cdot \nabla w)\} - sB(x)w.
$$
Hence by $u(\cdot,T) = 0$, we have
\begin{eqnarray*}
&& \int_Q \vert P_0u\vert^2 \weight dxdt
= \int_Q \vert Lw\vert^2 dxdt\\
=&& \int_Q \vert \ppp_tw + (H\cdot \nabla w)\vert^2 dxdt
+ \int_Q \vert sB\vert^2 w^2 dxdt
- 2s\int_Q Bw(\ppp_tw + (H\cdot\nabla w)) dxdt\\
\ge &&-2s \int_Q B(\ppp_tw + H\cdot \nabla w)w dxdt
+ s^2\int_Q B^2w^2 dxdt\\
= && -s \int_Q (B\ppp_t(w^2) + BH\cdot\nabla(w^2))dxdt
+ s^2\int_Q B^2w^2 dxdt\\
=&& s \int_Q  (\ppp_tB + (\mbox{div}\thinspace BH)) w^2 dxdt
- s\int^T_0\int_{\ppp\Omega} B(\nu\cdot H) w^2 dS_xdt\\
+ && s^2\int_Q B^2w^2 dxdt + s\int_{\o} B(x)\vert w(x,0)\vert^2 dx\\
\ge&& -M_0s\int_Q w^2 dxdt - s\int^T_0\int_{\ppp\Omega} B(\nu\cdot H) w^2
dS_xdt\\
+ && s^2\int_Q B^2w^2 dxdt + s\int_{\o} B(x)\vert w(x,0)\vert^2 dx.
\end{eqnarray*}
Substituting $w=e^{s\va}u$, we have
$$
s\int_{\OOO} B(x)\vert u(x,0)\vert^2 e^{2s\va(x,0)} dx
+ s^2\int_Q B^2(x)\vert u\vert^2 e^{2s\va} dxdt
$$
$$
\le \int_Q \vert P_0u\vert^2e^{2s\va} dxdt
+ M_0s\int_Q \vert u\vert^2 e^{2s\va} dxdt
+ s\int^T_0\int_{\ppp\OOO} B(x)\vert (H\cdot\nu)\vert \vert u\vert^2
e^{2s\va} dS_xdt.             \eqno{(2.5)}
$$
Next let $V \in L^{\infty}(\OOO)$, $\not\equiv 0$. Then
\begin{align*}
& \vert P_0u\vert^2 = \vert P_0u + Vu - Vu\vert^2
\le 2\vert P_0u+Vu\vert^2 + 2\vert Vu\vert^2\\
\le& 2\vert Pu\vert^2 + 2\Vert V\Vert^2_{L^{\infty}(\OOO)}\vert u\vert^2.
\end{align*}
Therefore (2.5) completes the proof of Lemma 1 (i).
\\
{\bf Proof of Lemma 1 (ii).}\\
By (2.2) and (2.3),
we have $B(x) \ge \mu - \beta > 0$ on $\ooo\o$.  Therefore we absorb
the second term on the right-hand side of the conclusion of Lemma 1 into the
left-hand side.
More precisely, we have
$$
s^2\int_Q B^2(x)\vert u\vert^2 e^{2s\va} dxdt
\ge (\mu-\beta)^2s^2\int_Q \vert u\vert^2 e^{2s\va} dxdt,
$$
and so we easily see that if $s\ge s_0$, then
\begin{align*}
& (\mu-\beta)^2s^2 - M_0s - 2\Vert V\Vert^2_{L^{\infty}(\OOO)}\\
\ge &\frac{(\mu-\beta)^2}{2}s^2
+ \left( \frac{(\mu-\beta)^2}{4}s^2 - M_0s\right)
+ \left( \frac{(\mu-\beta)^2}{4}s^2 - 2\Vert V\Vert^2_{L^{\infty}(\OOO)}
\right)\\
\ge & \frac{(\mu-\beta)^2}{2}s^2.
\end{align*}
Thus, by noting
$$
\int^T_0\int_{\ppp\o} B(x)(H\cdot\nu) \vert u\vert^2 dS_xdt
\le \int^T_0\int_{\ppp\o_+} B(x)(H\cdot\nu) \vert u\vert^2 dS_xdt,
$$
the proof of the part (ii) is completed.
\\
\vspace{0.2cm}

Next we show the classical energy estimate.
\\
{\bf Lemma 2.}\\
Let $\Vert V \Vert_{L^{\infty}(\OOO)}\le M$ and
$\Vert H\Vert_{\{C^1(\ooo\OOO)\}^n} \le M$ with arbitrarily fixed constant
$M>0$.  Let $w\in H^1(Q)$ satisfy
$$
\left\{ \begin{array}{rl}
&\ppp_tw + H(x)\cdot\nabla w + Vw = F(x,t) \quad
\mbox{in $Q$}, \\
& w(x,0) = a(x), \qquad x \in \OOO.\\
\end{array}\right.
\eqno{(2.6)}
$$
Then there exists a constant $C>0$, depending on $\Omega, T, M$, such that
$$
\int_{\Omega} \vert w(x,t)\vert^2 dx
+ \int^T_0 \int_{\ppp\OOO_+}(H\cdot\nu)\vert w\vert^2 dS_xdt
$$
$$
\le C\left( \Vert a\Vert^2_{L^2(\OOO)}
+ \Vert F\Vert^2_{L^2(Q)} + \int^T_0 \int_{\ppp\OOO_-}
\vert (H\cdot\nu)\vert \vert w\vert^2 dS_xdt\right), \quad 0\le t \le T.
                                                        \eqno{(2.7)}
$$
\\
{\bf Proof.}\\
Multiplying $\ppp_tw + H\cdot\nabla w + Vw = F$ by
$2w$ and integrating over $\OOO$, we have
$$
\ppp_t\int_{\OOO} \vert w(x,t)\vert^2 dx
+ \sum_{\ell=1}^n\int_{\OOO} h_{\ell}\ppp_{\ell}(\vert w\vert^2) dx
+ 2\int_{\OOO} V\vert w\vert^2 dx = 2\int_{\OOO} Fw dx.
$$
Setting $E(t) = \int_{\OOO} \vert w(x,t)\vert^2 dx$ and
integrating by parts, we obtain
\begin{eqnarray*}
&& E'(t) = - \int_{\ppp\OOO} (H\cdot \nu) \vert w\vert^2 dS_x
+ \int_{\OOO} (\ddd H) \vert w\vert^2 dx\\
-&& 2\int_{\OOO} Vw^2 dx + 2\int_{\OOO} Fw dx.
\end{eqnarray*}
Therefore, noting that $2\int_{\OOO} \vert Fw\vert dx \le
\int_{\OOO} \vert F\vert^2 dx + \int_{\OOO} \vert w\vert^2 dx$
and integrating over $(0,t)$, we have
$$
E(t) - E(0) = -\int^t_0\left( \int_{\ppp\o_+} + \int_{\ppp\o_-}\right)
(H\cdot \nu) \vert w\vert^2 dS_x dt
+ \int^t_0\int_{\OOO} (\mbox{div}\thinspace H)\vert w\vert^2 dxdt
$$
$$
- 2\int^t_0\int_{\OOO} Vw^2 dxdt + \int^t_0\int_{\OOO} \vert F\vert^2 dxdt
+ \int^t_0 E(\xi) d\xi.
$$
Therefore
$$
E(t) + \int^t_0\int_{\ppp\OOO_+} (H\cdot\nu) \vert w\vert^2 dS_xdt
\le \Vert a\Vert^2_{L^2(\OOO)} + \Vert F\Vert^2_{L^2(Q)}
$$
$$
+ (2\Vert V\Vert_{L^{\infty}(\OOO)} + \Vert \ddd H\Vert_{L^{\infty}
(\OOO)} + 1)\int^t_0 E(\xi) d\xi
+ \left\vert \int^t_0\int_{\ppp\OOO_-} (H\cdot\nu) \vert w\vert^2 dS_xdt
\right\vert                      \eqno{(2.8)}
$$
for $0 \le t \le T$.
Since $\int^t_0\int_{\ppp\OOO_+} (H\cdot\nu) \vert w\vert^2 dS_xdt
\ge 0$, $0\le t \le T$, omitting the second term on the left-hand side
of (2.8), we apply the Gronwall inequality, and we obtain
$$
E(t) \le Ce^{CT}\biggl(\Vert a\Vert^2_{L^2(\OOO)} + \Vert F\Vert^2_{L^2(Q)})
$$
$$
+ \int^T_0 \int_{\ppp\OOO_-} \vert (H\cdot\nu)\vert \vert w\vert^2
dS_xdt \biggr), \quad 0\le t \le T.                              \eqno{(2.9)}
$$
Substituting (2.9) into the third term on the right-hand side of (2.8),
we complete the proof of the lemma.
\section{Proofs of Theorems 1 and 2}~\\
{\bf Proof of Theorem 1.}\\
Henceforth $C>0$ denotes generic constants which are independent of
$s>0$.  We note that $Py = \ppp_ty + H\cdot \nabla y + Vy$ and
$\mu = \min_{x\in\ooo\OOO} (H(x)\cdot \nabla\psi(x)) > 0$, and we
set $R = \max_{x\in\ooo\o} \psi(x)$ and
$r = \min_{x\in\ooo\o} \psi(x)$.  By (1.7), we can choose
$\beta > 0$ such that
$$
T > \frac{R-r}{\beta}, \quad 0 < \beta < \mu.    \eqno{(3.1)}
$$
In fact, it suffices to choose $\beta > 0$ such that $\beta\in (0,\mu)$
is sufficiently close to $\mu$.

With this $\beta>0$, we set
$$
\va(x,t) = -\beta t + \psi(x), \qquad (x,t) \in Q.    \eqno{(3.2)}
$$
Then (3.1) implies
$$
\va(x,T) \le R - \beta T < r \le \va(x',0), \quad x, x' \in \ooo\o.
$$
Therefore by $\va\in C^1(\ooo{Q)}$, there exist $\delta_1>0$ and $r_0, r_1$
such that $R - \beta T < r_0 < r_1 < r$,
$$
\left\{ \begin{array}{rl}
&\va(x,t) > r_1, \quad x \in \ooo\o, \thinspace 0 \le t\le \delta_1,\\
&\va(x,t) < r_0, \quad x \in \ooo\o, \thinspace T-2\delta_1 \le t\le T.
\end{array}\right.
                                            \eqno{(3.3)}
$$
For applying Lemma 1, we need a cut-off function $\chi\in C^{\infty}_0
({\Bbb R})$ such that $0 \le \chi \le 1$ and
$$
\chi(t)  =
\left\{
\begin{array}{rl}
1, \qquad & 0 \le t \le T-2\delta_1,\\
0, \qquad & T-\delta_1\le t \le T.
\end{array}\right.
\eqno{(3.4)}
$$
We set
$$
z = (\ppp_ty)\chi.
$$
Then $z(x,T) = 0$, $x \in \OOO$ and
$$
Pz = \chi f(\ppp_tR) + (\ppp_t\chi)\ppp_ty, \quad (x,t) \in Q
$$
and
$$
z(x,0) = f(x)R(x,0), \qquad x\in\o.
$$
Applying Lemma 1 (ii) to $z$, we obtain
$$
s\int_{\OOO} \vert z(x,0)\vert^2 e^{2s\va(x,0)} dx
\le C\int_Q \vert \chi f(\ppp_tR)\vert^2 \weight dxdt
+ C\int_Q \vert (\ppp_t\chi)\ppp_ty\vert^2 \weight dxdt
+ Ce^{Cs}D^2                        \eqno{(3.5)}
$$
for all large $s > 0$.
Here
$$
D^2 = \int^T_0\int_{\ppp\Omega_+} (H\cdot\nu)\vert \ppp_ty\vert^2 dS_xdt.
$$
Since $\ppp_t\chi = 0$ for $0 \le t\le T-2\delta_1$ or $T-\delta_1
\le t \le T$, by (3.3) and the a priori bound $\Vert \ppp_ty\Vert_{L^2(Q)}
\le M$, we have
$$
\int_Q \vert (\ppp_t\chi)\ppp_ty\vert^2 \weight dxdt
= \int^{T-\delta_1}_{T-2\delta_1} \int_{\OOO}
\vert (\ppp_t\chi)\ppp_ty\vert^2 \weight dxdt
\le Ce^{2sr_0}\int^{T-\delta_1}_{T-2\delta_1} \int_{\OOO}
\vert \ppp_ty\vert^2 dxdt                     \eqno{(3.6)}
$$
and
$$
\int_Q \vert (\ppp_t\chi)\ppp_ty\vert^2 \weight dxdt
\le Ce^{2sr_0}M^2                                       \eqno{(3.7)}
$$
for all large $s>0$.
Moreover $R(x,0) \ne 0$ for $x \in \ooo\OOO$ and
$z(x,0) = f(x)R(x,0)$, $x \in \ooo\OOO$, we have
$$
\int_{\OOO} \vert z(x,0)\vert^2 e^{2s\va(x,0)} dx
\ge C\int_{\OOO} \vert f(x)\vert^2 e^{2s\va(x,0)} dx.
$$
Therefore (3.5) yields
$$
s\int_{\OOO} \vert f(x)\vert^2 e^{2s\va(x,0)} dx
\le C\int_Q \vert f(x)\vert^2 e^{2s\va(x,t)} dxdt
+ CM^2e^{2sr_0} + Ce^{Cs}D^2.
$$
Since $\va(x,t) \le \va(x,0)$ for $(x,t) \in Q$, we have
\begin{align*}
&s\int_{\OOO} \vert f(x)\vert^2 e^{2s\va(x,0)} dx
\le C\int^T_0\int_{\OOO} \vert f(x)\vert^2 e^{2s\va(x,0)} dxdt
+ CM^2e^{2sr_0} + Ce^{Cs}D^2\\
\le &CT\int_{\OOO} \vert f(x)\vert^2 e^{2s\va(x,0)} dx
+ CM^2e^{2sr_0} + Ce^{Cs}D^2,
\end{align*}
that is,
$$
(s-CT)\int_{\OOO} \vert f(x)\vert^2 e^{2s\va(x,0)} dx
\le CM^2e^{2sr_0} + Ce^{Cs}D^2
$$
for all large $s>0$.  Using $\va(x,0) > r_1$ by (3.3)
and choosing $s>0$ large,
we obtain
$$
se^{2sr_1}\int_{\OOO} \vert f(x)\vert^2 dx
\le CM^2e^{2sr_0} + Ce^{Cs}D^2,
$$
that is,
$$
\Vert f\Vert^2_{L^2(\OOO)}
\le CM^2e^{-2sr_*} + Ce^{Cs}D^2          \eqno{(3.8)}
$$
for all large $s > s_*$, where $s_*>0$ is a sufficiently large
constant.  Here we set $r_* := r_1 - r_0 > 0$.
We separately consider the two cases: $D \ge M$ and $D<M$.
\\
{\bf Case 1 $D\ge M$:} \\
Estimate (3.8) implies
$$
\Vert f\Vert^2_{L^2(\OOO)} \le (Ce^{-2sr_*} + Ce^{Cs})D^2.       \eqno{(3.9)}
$$
\\
{\bf Case 2 $D<M$:}\\
Replacing $C$ by $Ce^{Cs_*}$, we see that (3.8) holds for all $s>0$.
We make the right-hand side of (3.8) small in $s$.  We choose
$M^2e^{-2sr_*} = e^{Cs}D^2$, that is,
$$
s = \frac{2}{C+2r_*}\log \frac{M}{D}.
$$
Therefore (3.8) is reduced to
$$
\Vert f\Vert_{L^2(\OOO)}\le 2CM^{1-\theta}D^{\theta},
$$
where $\theta = \frac{2r_*}{C+2r_*} \in (0,1)$.
Thus the proof of Theorem 1 (i) is completed.

Next we prove Theorem 1 (ii).  We have
$$
\left\{\begin{array}{rl}
& \ppp_t(\ppp_ty) + H(x)\cdot\nabla \ppp_ty + V\ppp_ty = f(x)\ppp_tR,
\quad (x,t) \in Q,\\
& (\ppp_ty)(x,0) = f(x)R(x,0), \qquad x\in \Omega.
\end{array}\right.
$$
Applying Lemma 2 to $\ppp_ty$, we obtain
$$
\int_{\OOO} \vert \ppp_ty(x,t)\vert^2 dx
+ \int^T_0\int_{\ppp\OOO_+} (H\cdot\nu)\vert \ppp_ty\vert^2 dS_xdt
$$
$$
\le C\Vert f\Vert^2_{L^2(\OOO)}
+ C\int^T_0\int_{\ppp\OOO_-} \vert (H\cdot\nu)\vert
\vert \ppp_ty\vert^2 dS_xdt                      \eqno{(3.10)}
$$
for $0 \le t \le T$.  Therefore, omitting the second term on the
left-hand side of (3.10) and applying it to (3.6), we obtain
$$
\int_Q \vert (\ppp_t\chi)\ppp_ty\vert^2 \weight dxdt
\le Ce^{2sr_0}\Vert f\Vert^2_{L^2(\OOO)}
+ Ce^{2sr_0}\int^T_0\int_{\ppp\OOO_-}
\vert (H\cdot\nu)\vert\vert \ppp_ty\vert^2 dS_xdt
$$
and similarly to (3.8), from (3.5) we can obtain
$$
\Vert f\Vert^2_{L^2(\OOO)} \le Ce^{-2sr_*}\Vert f\Vert^2_{L^2(\OOO)}
+ Ce^{Cs}\int^T_0\int_{\ppp\OOO} \vert (H\cdot\nu)\vert
\vert \ppp_ty\vert^2 dS_xdt.
$$
Choosing $s>0$ large, we can absorb the first term on the right-hand side
into the left-hand side, and complete the proof of (ii).

Finally we prove (iii).  By (1.8), the conclusion of (ii) immediately
yields the second inequality in (1.9).  Next, by (1.8) and (3.10), we have
$$
\int^T_0\int_{\ppp\OOO_+}
\vert (H\cdot\nu)\vert\vert \ppp_ty\vert^2 dS_xdt
\le C\Vert f\Vert^2_{L^2(\OOO)},
$$
which proves the first inequality in (1.9).  Thus the proof of
Theorem 1 (iii) is completed.
\\
\vspace{0.2cm}

{\bf Proof of Theorem 2}.\\
Theorem 2 can be derived directly by Theorem 1.  In fact,
setting $y = u_1 - u_2$, $f = V_1-V_2$ and $R=-u_2$, by
(1.10) - (1.14) we have
$$
\left\{\begin{array}{rl}
& \ppp_ty + H(x)\cdot\nabla y + V_1y = f(x)R(x,t),
\quad (x,t) \in Q,\\
& y(x,0) = 0, \qquad x\in \Omega,\\
& y\vert_{\ppp\Omega_-\times (0,T)} = 0,
\end{array}\right.
$$
and $\ppp_ty \in H^1(Q)$, $\Vert \ppp_tR\Vert_{L^2(0,T;L^{\infty}
(\OOO))} \le M$, $R(x,0) = -a(x) \ne 0$ for $x \in \ooo\OOO$.
Thus Theorem 1 (iii) yields the conclusion of Theorem 2, and the
proof of Theorem 2 is completed.
\section{Proof of Theorem 3}~\\
We set
$$
P_du = \ppp_tu + \nabla d \cdot \nabla u + u\Delta d.
$$
For the proof, we further need a Carleman estimate for $\Delta$ with the same
weight function for $P_d$.  Unfortunately the weight function defined
by (2.1) does not work as weight for a Carleman estimate for $\Delta$.
Thus we have to introduce a second large parameter in the weight
function.  That is, as the weight function, we set
$$
{\va}_d(x,t) = e^{\lambda(-\beta t + d(x))}, \quad (x,t) \in Q,
                                              \eqno{(4.1)}
$$
where $\la>0$ is chosen later.  The weight function in the form (4.1)
has been known as more flexible weight function producing
a Carleman estimate (e.g., H\"ormander [10], Isakov [13]).

In the existing works on inverse problems by Carleman estimates,
given weight functions have been used, and to the
best knowledge of the authors,  Theorem 3 is the first case where a Carleman
estimate with unknown coefficient as the weight, is seriously involved.
The use of such a Carleman estimate is necessary in order that the
admissible set $\mathcal{D}$ of unknown coefficients $d$'s is generously
formulated such as (1.17).
Otherwise the admissible set of unknown coefficients should be more
restrictive.  For example, for inverse problems of determining
principal parts for second-order hyperbolic equations, there are no
works by Carleman estimate with weight function given by unknown
coefficient, so that choices of admissible sets of unknown coefficients of the
principal parts are very limited (e.g., Bellassoued and Yamamoto
[6], Section 6 of Chapter 5).
For it, we need a Carleman estimate where all the constants can be taken
uniformly for arbitrary $d \in \mathcal{D}$.

We first show a Carleman estimate for $P_d$ with weight (4.1), which should
hold uniformly for all $d \in \mathcal{D}$.
We fix $\beta, T>0$ such that
$$
0 < \beta < \delta_0^2, \quad T > \frac{\sup_{d\in \mathcal{D}}
(\max_{x\in\ooo\OOO} d(x) - \min_{x\in\ooo\OOO} d(x))}{\beta},    \eqno{(4.2)}
$$
where the constant $\delta_0$ characterizes $\mathcal{D}$ (see (1.17)).
Henceforth we set
$$
J_d(x,t) := \ppp_t{\va}_d + (\nabla d\cdot \nabla {\va}_d)
= \la{\va}_d(-\beta + \vert \nabla d\vert^2).
$$
\\
{\bf Lemma 3}\\
For each $\la>0$ there exists a constant $s_0 > 0$, which is
dependent on $\la, \delta_0, M$, such that
we can choose a constant $C>0$ satisfying
$$
\int_{\OOO} s\la\va_{d}(x,0)\vert u(x,0)\vert^2 e^{2s\va_{d}(x,0)} dx
+ \int_Q s^2\la^2\va_{d}^2 \vert u(x,t)\vert^2 e^{2s\va_{d}} dxdt
$$
$$
\le C\int_Q \vert P_du\vert^2 e^{2s\va} dxdt
+ \int^T_0\int_{\ppp\OOO} sJ_d(\ppp_{\nu}d) \vert u\vert^2 e^{2s\va_{d}} dS_xdt
                                        \eqno{(4.3)}
$$
for all $s \ge s_0$, all $d\in \mathcal{D}$ and all
$u \in H^1(Q)$ satisfying $u(x,T) = 0$, $x \in \OOO$.
\\
\vspace{0.3cm}

Here we note that the choices of $C$ and $s_0$ in (4.3) are uniformly for
$d \in \mathcal{D}$.
\\
{\bf Proof}.\\
The proof is similar to Lemma 1 by noting the independency of the
constant of $d \in \mathcal{D}$.  Let $d \in \mathcal{D}$ be chosen
arbitrarily.  First we consider the case of
$P_du = \ppp_tu + \nabla d \cdot \nabla u$.  We set
$w = e^{s\va_{d}}u$, and
$$
Lw = e^{s\va_{d}}P_d(e^{-s\va_{d}}w)
= \ppp_tw + \nabla d\cdot \nabla w - sJ_d(x,t)w.
$$
We note $\ppp_tJ_d(x,t) = -\la^2\beta{\va}_d(-\beta + \vert \nabla d
\vert^2)$,
$\ppp_kJ_d(x,t) = \la^2\va_{d}(\ppp_kd)(-\beta + \vert \nabla d\vert^2)
+ \la\va_{d}(\ppp_k\vert \nabla d \vert^2)$ and
$$
J_d(x,t) \ge \la\va_{d}(\delta_0^2-\beta), \quad (x,t) \in Q, \thinspace
d \in \mathcal{D}                             \eqno{(4.4)}
$$
by (1.17).
Then
\begin{align*}
& \int_Q \vert P_du\vert^2 e^{2s\va_{d}} dxdt
= \int_Q \vert Lw\vert^2 dxdt\\
=& \int_Q \vert (\ppp_tw + \nabla d\cdot\nabla w) - sJ_dw\vert^2 dxdt\\
=& \int_Q \vert \ppp_tw + \nabla d\cdot\nabla w\vert^2 dxdt
+ s^2\int_Q J_d^2\vert w\vert^2 dxdt
- 2s\int_Q (J_dw)(\ppp_tw + \nabla d\cdot\nabla w) dxdt\\
\ge& s^2\int_Q J_d^2\vert w\vert^2 dxdt
- s\int_Q J_d(2w(\ppp_tw)) dxdt
- s\int_Q J_d\nabla d \cdot (2w\nabla w) dxdt\\
=: & I_1 + I_2 + I_3.
\end{align*}

Henceforth $C>0$ denotes generic constants which are independent of
$d \in \mathcal{D}$ and $s>0$, but dependent on $\delta_0, M, \la$.
Next integration by parts yields
\begin{align*}
& I_2 = -s\int_Q J_d\ppp_t(\vert w\vert^2) dxdt
= -s\int_{\OOO} [J_d\vert w\vert^2]^T_0 dx
+ s\int_Q (\ppp_tJ_d)\vert w\vert^2 dxdt \\
\ge& s\int_{\OOO} \la\va_{d}(x,0)(-\beta + \vert\nabla d\vert^2)
\vert w(x,0)\vert^2 dx
- C\int_Q s\la^2\va_{d}\vert w\vert^2 dxdt\\
\ge& s(\delta_0^2-\beta)\int_{\OOO} s\la\va_{d}(x,0)\vert w(x,0)\vert^2 dx
- C\int_Q s\la^2\va_{d}\vert w\vert^2 dxdt
\end{align*}
and
\begin{align*}
& I_3 = -s\int_Q \sum_{k=1}^n (J_d\ppp_kd)\ppp_k(\vert w\vert^2) dxdt\\
=& -s\int^T_0\int_{\ppp\OOO} \sum_{k=1}^n J_d(\ppp_kd)\nu_k\vert w\vert^2
dS_xdt
+ s\int_Q \sum_{k=1}^n \ppp_k(J_d\ppp_kd)\vert w\vert^2 dxdt\\
\ge& -\int^T_0 \int_{\ppp\OOO} sJ_d(\ppp_{\nu}d)\vert w\vert^2 dS_xdt
- C\int_Q s\la^2\va_{d} \vert w\vert^2 dxdt.
\end{align*}
Hence
$$
\int_{\OOO} s\la\va_{d}(x,0)\vert w(x,0)\vert^2 dx
+ \int_Q s^2\la^2\va_{d}^2\vert w\vert^2 dxdt
$$
$$
\le C\int_Q \vert Lw\vert^2 dxdt + C\int_Q s\la^2\va_{d} \vert w\vert^2 dxdt
+ \int^T_0\int_{\ppp\OOO} sJ_d(\ppp_{\nu}d)\vert w\vert^2 dS_xdt.  \eqno{(4.5)}
$$
We choose
$$
s_0 > 2Ce^{\la(\beta T + M)}.
$$
Then, since $\va_{d} (x,t) \ge e^{-\la\beta T - \la M}$ for
$(x,t) \in Q$ and $d \in \mathcal{D}$, if $s \ge s_0$, then
\begin{align*}
& s^2\la^2 \va_{d}^2 - Cs\la^2\va_{d}
= s^2\la^2\va_{d}^2\left( 1 - \frac{C}{s\va_{d}}\right)\\
\ge & s^2\la^2\va_{d}^2\left( 1 - \frac{Ce^{\la(\beta T+M)}}{s_0}\right)
\ge \frac{1}{2}s^2\la^2\va_{d}^2.
\end{align*}
Therefore (4.5) yields
$$
\int_{\OOO} s\la\va_{d}(x,0)\vert w(x,0)\vert^2 dx
+ \frac{1}{2}\int_Q s^2\la^2\va_{d}^2\vert w\vert^2 dxdt
$$
$$
\le  C\int_Q \vert Lw\vert^2 dxdt
+ \int^T_0\int_{\ppp\OOO} sJ_d(\ppp_{\nu}d)\vert w\vert^2 dS_xdt
                                                             \eqno{(4.6)}
$$
for all $s \ge s_0$.  Henceforth $s_0>0$ denotes
generic constants which are dependent on $\la, \delta_0, M$.
Since $w = ue^{s\va_{d}}$, we rewrite (4.6) in view of $u$, so that
$$
\int_{\OOO} s\la\va_{d}(x,0)\vert u(x,0)\vert^2 e^{2s\va_{d}(x,0)}dx
+ \int_Q s^2\la^2\va_{d}^2\vert u\vert^2 e^{2s\va_{d}} dxdt
$$
$$
\le C\int_Q \vert P_du\vert^2 e^{2s\va_{d}} dxdt
+ \int^T_0\int_{\ppp\OOO} sJ_d(\ppp_{\nu}d)\vert u\vert^2e^{2s\va_{d}} dS_xdt
$$
for all $s \ge s_0$.  Here we note that $C>0$ and $s_0>0$ may vary
line by line, but dependent only on $\la, \delta_0, M$.
Next
\begin{align*}
& \vert \ppp_tu + \nabla d\cdot\nabla u + (\Delta d)u\vert^2
\le 2\vert \ppp_tu + \nabla d\cdot\nabla u\vert^2
+ 2\vert (\Delta d)u\vert^2\\
\le & 2\vert \ppp_tu + \nabla d\cdot\nabla u\vert^2
+ 2M^2\vert u\vert^2
\end{align*}
in $Q$ for each $d \in \mathcal{D}$.  Similarly to the proof of Lemma 1 (i),
we can finish the proof of (4.3).
\\
\vspace{0.3cm}

Next we show a Carleman estimate for $\Delta$ with the weight function
$\va_{d}(x,0)$.
\\
{\bf Lemma 4}\\
There exists a constant $\la_0 = \la_0(\delta_0,M) > 0$ such that
for any $\la \ge \la_0$, there exists a constant $s_1
= s_1(\la,\delta_0,M)>0$ satisfying: we can choose $C>0$ such that
\begin{align*}
&\int_{\OOO} (s\la^2\va_{d}(x,0)\vert \nabla f(x)\vert^2
+ s^3\la^4\va_{d}(x,0)^3\vert f(x)\vert^2)e^{2s\va_{d}(x,0)} dx\\
\le& C\int_{\OOO} \vert \Delta f\vert^2 e^{2s\va_{d}(x,0)} dx
\end{align*}
for all $s \ge s_1$, $d\in \mathcal{D}, f \in H^2_0(\OOO)$.
\\
\vspace{0.2cm}

Lemma 4 is a Carleman estimate for $\Delta$ and the Carleman
estimate for $\Delta$ is well known
(e.g., H\"ormander [10], Isakov [13]).  However we need the uniformity
of the constants $s_0, C$ with respect to $d \in \mathcal{D}$
in the Carleman estimate.  Thus in Appendix we
give a direct proof by integration by parts, which differs from
[10], [13].
\\
\vspace{0.2cm}

Now we proceed to \\
{\bf Proof of Theorem 3}.\\
{\bf First Step.} We set
$$
m_0 = \sup_{d\in \mathcal{D}} (\max_{x\in \ooo\OOO} d(x)
- \min_{x\in \ooo\OOO} d(x)).
$$
By assumption (1.21): $T > \frac{m_0}{\delta_0^2}$ of the theorem, for
arbitrary $d\in \mathcal{D}$, we can choose $\beta > 0$ and fix such that
$$
0 < \beta < \delta_0^2, \qquad T > \frac{m_0}{\beta}.   \eqno{(4.7)}
$$
We choose $\la_0 > 0$ given in Lemma 4 and fix.
Next we set $s_* = \max\{ s_0(\la_0,\delta_0,M),
s_1(\la_0,\delta_0,M)\}$, where $s_j(\la_0,\delta_0,M)$, $j=0,1$ are
proved to exist in Lemmata 3 and 4.
Since $\ppp_{\nu}d = g_2$ on $\ppp\OOO$ for any $d\in
\mathcal{D}$, we note that the subboundaries $\Gamma_+$ and $\Gamma_-$
defined by (1.18) correspond to $\ppp\Omega_+$ and $\ppp\Omega_-$ respectively.
Since $\la_0>0$ is fixed,
henceforth by $C$ denoting generic constants which are independent of
$s \ge s_*$ and $d\in \mathcal{D}$ but dependent on $\la_0$, $a$,
$\beta$, $m_0$, $\Omega$, $T$, noting (4.7) and (4.4),
by Lemmata 3 and 4 we obtain
$$
\int_{\OOO} s\vert y(x,0)\vert^2 e^{2s\va_{d}(x,0)} dx
+ \int_Q s^2\vert y\vert^2 e^{2s\va_{d}} dxdt
$$
$$
\le  C\int_Q \vert P_dy\vert^2 e^{2s\va_{d}} dxdt
+ Ce^{Cs}\int^T_0\int_{\Gamma_+} \vert y\vert^2 dS_xdt  \eqno{(4.8)}
$$
and
$$
\int_{\OOO} (s\vert \nabla f\vert^2 + s^3\vert f\vert^2)e^{2s\va_{d}(x,0)} dx
\le C\int_{\OOO} \vert \Delta f\vert^2 e^{2s\va_{d}(x,0)} dx  \eqno{(4.9)}
$$
for all $s \ge s_*$, $d\in \mathcal{D}$, $y \in H^1(Q)$ satisfying
$y(\cdot,T)=0$ in $\OOO$ and $f \in H^2_0(\OOO)$.
Thus we obtained two Carleman estimates which hold uniformly for
arbitrary $d \in \mathcal{D}$.
\\
{\bf Second Step.}  We set
$y = \rho_1 - \rho_2$, $f = d_1-d_2$ and $R(x,t) = -\rho_2(x,t)$ in $Q$.
Then, by (1.19) and (1.20), we obtain
$$
\left\{ \begin{array}{rl}
& \ppp_ty + \nabla d_1\cdot \nabla y + y\Delta d_1
= \nabla f\cdot \nabla R + R\Delta f \quad \mbox{in $Q$},\\
& y(x,0) = 0, \qquad x\in \OOO,\\
& y=0 \qquad \mbox{on $\Gamma_-\times (0,T)$}.
\end{array}\right.
$$
We differentiate this partial differential equation with respect to $y$
in $t$, noting that
$\ppp_ty \in H^1(Q)$, $\ppp_tR \in L^2(0,T;W^{1,\infty}(\OOO))$
and $R(x,0) = -a(x)$.  Then we have
$$
\left\{ \begin{array}{rl}
& \ppp_t(\ppp_ty) + \nabla d_1\cdot \nabla (\ppp_ty) + (\ppp_ty)\Delta d_1\\
= & \nabla f\cdot \nabla\ppp_t R + (\ppp_tR)\Delta f \quad \mbox{in $Q$},\\
& \ppp_ty(x,0) = -\nabla f\cdot \nabla a - a\Delta f, \quad x\in \OOO,\\
& \ppp_ty=0 \qquad \mbox{on $\Gamma_-\times (0,T)$}.
\end{array}\right.
\eqno{(4.10)}
$$
Applying Lemma 2 to (4.10) and estimating $\ppp_ty$, we readily verify
the first inequality in (1.23).
Thus the rest part of this section is devoted to the proof of
the second inequality of (1.23).

Since $y(\cdot,T) = 0$ does not hold, for applying (4.8), we need
a cut-off function.  Noting
$$
\vadd(x,0) = e^{\la_0d_1(x)} \ge e^{\la_0 \min_{x\in\ooo\OOO} d_1(x)}
$$
and
$$
\vadd(x,T) = e^{\la_0(-\beta T + d_1(x))} \le e^{-\la_0\beta T}
e^{\la_0 \max_{x\in\ooo\OOO} d_1(x)}
$$
for all $x \in \ooo\OOO$.
Henceforth we set
$$
r_0 = \frac{1}{2} \la_0e^{-\la_0 \beta - \la_0M}(\beta T - m_0).
$$
By (4.7) we have $r_0 > 0$. We note that $r_0>0$ is independent of $s>0$ and
special choices of
$d_1, d_2\in \mathcal{D}$, but dependent on $\delta_0, M, \la_0$.

The mean value theorem implies
\begin{align*}
& \max_{x\in\ooo\OOO} \vadd(x,T) - \min_{x\in\ooo\OOO} \vadd(x,0)
= e^{-\la_0\beta T + \la_0\max_{x\in\ooo\OOO} d_1(x)} -
e^{\la_0 \min_{x\in\ooo\OOO} d_1(x)}\\
=& e^{\xi_1}\la_0 \{ (\max_{x\in\ooo\OOO} d_1(x)
- \min_{x\in\ooo\OOO} d_1(x)) - \beta T\}
\le \la_0e^{\xi_1}(m_0 - \beta T),
\end{align*}
where $\xi_1$ is a number between
$\la_0\min_{x\in\ooo\OOO} d_1(x)$ and $-\la_0\beta T
+ \la_0\max_{x\in\ooo\OOO} d_1(x)$.

Therefore, by $d \in \mathcal{D}$, we have
$e^{\xi_1} \ge e^{-\la_0\beta T - \la_0M}$, and
$$
\max_{x\in\ooo\OOO} \vadd(x,T) - \min_{x\in\ooo\OOO} \vadd(x,0)
\le -\la_0e^{\xi_1}(\beta T - m_0) \le -\la_0e^{-\la_0\beta T-\la_0M}
(\beta T - m_0) = -2r_0.                                  \eqno{(4.11)}
$$
For any $d_1 \in \mathcal{D}$, we have
$$
\vert \max_{x\in\ooo\OOO} \va_{d_1}(x,T)
- \max_{x\in\ooo\OOO} \va_{d_1}(x,t)\vert
= \vert \exp(\la_0(-\beta T + \max_{x\in\ooo\OOO} d_1(x)))
- \exp(\la_0(-\beta t + \max_{x\in\ooo\OOO} d_1(x))\vert
$$
$$
= \vert \exp(\la_0\max_{x\in\ooo\OOO} d_1(x))\vert
\vert e^{-\la_0\beta T} - e^{-\la_0\beta t}\vert
\le e^{\la_0 M}\la_0\beta \vert T-t\vert.                    \eqno{(4.12)}
$$
At the last inequality, by $d_1 \in \mathcal{D}$ and the mean value
theorem, we have
$$
\max_{x\in\ooo\OOO} d_1(x) \le M
$$
and
$$
\vert e^{-\la_0\beta T} - e^{-\la_0\beta t}\vert
= \vert e^{-\xi_2}(-\la_0\beta T + \la_0\beta t)\vert
\le \la_0\beta \vert T-t\vert,
$$
where $\xi_2\ge 0$ is some constant.

We fix a constant $\delta_1 > 0$ sufficienly small such that
$$
2e^{\la_0M}\la_0\beta \delta_1 < r_0.
$$
Then, for $T-2\delta_1 \le t \le T$ and any $d_1 \in \mathcal{D}$,
by (4.11) and (4.12), we have
\begin{align*}
& \max_{x\in\ooo\OOO} \va_{d_1}(x,t) - \min_{x\in\OOO} \va_{d_1}(x,0)\\
=& \max_{x\in\ooo\OOO} \va_{d_1}(x,t) - \max_{x\in\ooo\OOO} \va_{d_1}(x,T)
+ \max_{x\in\ooo\OOO} \va_{d_1}(x,T) - \min_{x\in\ooo\OOO} \va_{d_1}(x,0)\\
\le& \vert \max_{x\in\ooo\OOO} \va_{d_1}(x,t)
- \max_{x\in\ooo\OOO} \va_{d_1}(x,T)\vert - 2r_0
< e^{\la_0M}\la_0\beta \vert t-T\vert - 2r_0\\
\le& 2e^{\la_0M}\la_0 \beta \delta_1 - 2r_0
< r_0 -2r_0 = -r_0.
\end{align*}
Therefore
$$
\max_{x\in\ooo\OOO} \va_{d_1}(x,t) < \mu_0 - r_0,
\quad T-2\delta_1 \le t \le T,                         \eqno{(4.13)}
$$
where we set $\mu_0 : = \sup_{d\in \mathcal{D}}\min_{x\in\OOO} e^{\la_0d(x)}$.

Let $\chi \in C^{\infty}_0({\Bbb R})$ satisfy $0 \le \chi \le 1$ and (3.4)
with here chosen $\delta_1$.
We set
$$
z = \chi\ppp_ty.
$$
Then, by (4.10), we have
$$
\left\{ \begin{array}{rl}
& P_{d_1}z = \ppp_tz + \nabla d_1\cdot \nabla z + z\Delta d_1\\
= & \chi\nabla f\cdot \nabla\ppp_t R + \chi(\ppp_tR)\Delta f
+ (\ppp_t\chi)\ppp_ty \quad \mbox{in $Q$},\\
& z(x,0) = -\nabla f\cdot \nabla a - a\Delta f, \quad z(x,T) = 0
\quad x\in \OOO,\\
& z = 0 \qquad \mbox{on $\Gamma_-\times (0,T)$}.
\end{array}\right.
\eqno{(4.14)}
$$
Applying (4.8) to (4.14), we obtain
$$
\int_{\OOO} s\vert a\Delta f + \nabla f\cdot\nabla a\vert^2
e^{2s\vadd(x,0)} dx
$$
$$
\le C\int_Q \vert\chi\nabla f\cdot \nabla\ppp_t R + \chi(\ppp_tR)\Delta f
\vert^2 e^{2s\vadd} dxdt
+ C\int_Q \vert (\ppp_t\chi)\ppp_ty\vert^2 e^{2s\vadd} dxdt
+ Ce^{Cs}D^2                \eqno{(4.15)}
$$
for all $s \ge s_*$.  Here and henceforth we set
$$
D^2 = \int^T_0 \int_{\Gamma_+} \vert \ppp_t(\rho_1-\rho_2)\vert^2
dS_xdt.
$$
Apply Lemma 2 to (4.10), by (1.21) we have
$$
\int_{\OOO} \vert \ppp_ty(x,t)\vert^2 dx
\le C(\Vert \nabla f\cdot \nabla a + a\Delta f\Vert^2_{L^2(\OOO)}
+ \Vert \nabla f\cdot \nabla (\ppp_tR) + (\ppp_tR)\Delta f\Vert^2_{L^2(Q)})
+ CD^2
$$
$$
\le C(\Vert \Delta f\Vert^2_{L^2(\OOO)} + \Vert \nabla f\Vert^2
_{\{L^2(\OOO)\}^n} + D^2), \quad 0\le t \le T.            \eqno{(4.16)}
$$
Consequently, since $\ppp_t\chi \ne 0$ only if $T-2\delta_1 \le t \le
T-\delta_1$ by (3.4), it follows from (1.14) that (4.15) and (4.16)
yield
\begin{align*}
& \int_{\OOO} s\vert \Delta f\vert^2 e^{2s\vadd(x,0)} dx
- C\int_{\OOO} s\vert \nabla f\vert^2 e^{2s\vadd(x,0)} dx\\
\le & C\int_{\OOO} (\vert \Delta f\vert^2 + \vert \nabla f\vert^2)
e^{2s\vadd(x,0)} dx\\
+& C\left( \int^{T-\delta_1}_{T-2\delta_1}\int_{\OOO}
e^{2s\vadd} dxdt\right) (\Vert \Delta f\Vert^2_{L^2(\OOO)}
+ \Vert \nabla f\Vert^2_{\{L^2(\OOO)\}^n} + D^2)
+ Ce^{Cs}D^2
\end{align*}
for all $s \ge s_*$.
For the first integral on the right-hand side, we used
$\vadd(x,t) \le \vadd(x,0)$ for $(x,t) \in Q$, and so
\begin{align*}
&\int_Q \vert\chi\nabla f\cdot \nabla(\ppp_t R) + \chi(\ppp_tR)\Delta f
\vert^2 e^{2s\vadd(x,t)} dxdt \\
\le& C\int_Q (\vert \nabla f\vert^2 + \vert \Delta f\vert^2)
e^{2s\vadd(x,t)} dxdt
\le C\int_{\OOO} (\vert \nabla f\vert^2 + \vert \Delta f\vert^2)
e^{2s\vadd(x,0)} dx.
\end{align*}

Therefore we obtain
\begin{align*}
& \int_{\OOO} s\vert \Delta f\vert^2 e^{2s\vadd(x,0)} dx\\
\le &C\int_{\OOO} \vert \Delta f\vert^2 e^{2s\vadd(x,0)} dx
+ Cs\int_{\OOO} \vert \nabla f\vert^2 e^{2s\vadd(x,0)} dx\\
+ & C\max_{T-2\delta_1\le t\le T-\delta_1}
\max_{x\in \ooo\OOO} e^{2s\vadd(x,t)}(\Vert \Delta f\Vert^2_{L^2(\OOO)}
+ \Vert \nabla f\Vert^2_{\{L^2(\OOO)\}^n})\\
+& Ce^{Cs}D^2
\end{align*}
for all $s \ge s_*$.  The first term on the right-hand side can be
absorbed into the left-hand side by choosing $s>0$ large,
and by (4.13) we have
$$
\int_{\OOO} s\vert \Delta f\vert^2 e^{2s\vadd(x,0)} dx
$$
$$
\le Cs\int_{\OOO} \vert \nabla f\vert^2 e^{2s\vadd(x,0)} dx
+  Ce^{2s(\mu_0-r_0)}(\Vert \Delta f\Vert^2_{L^2(\OOO)}
+ \Vert \nabla f\Vert^2_{\{L^2(\OOO)\}^n}) + Ce^{Cs}D^2     \eqno{(4.17)}
$$
for all $s \ge s_*$.

Since $f = d_1 - d_2 \in H^2_0(\OOO)$ by $d_1, d_2 \in
\mathcal{D}$, we can apply (4.9) to obtain
$$
\int_{\OOO} (s^2\vert \nabla f\vert^2 + s^4\vert f\vert^2)
e^{2s\vadd(x,0)} dx
\le C\int_{\OOO} s\vert \Delta f\vert^2 e^{2s\vadd(x,0)} dx.
$$
Substituting this into the left-hand side of (4.17), we obtain
\begin{align*}
& \int_{\OOO} (s\vert \Delta f\vert^2 + s^2\vert \nabla f\vert^2
+ s^4\vert f\vert^2)e^{2s\vadd(x,0)} dx\\
\le &Cs\int_{\OOO} \vert \nabla f\vert^2 e^{2s\vadd(x,0)} dx
+  Ce^{2s(\mu_0-r_0)}(\Vert \Delta f\Vert^2_{L^2(\OOO)}
+ \Vert \nabla f\Vert^2_{\{L^2(\OOO)\}^n}) + Ce^{Cs}D^2
\end{align*}
for all $s \ge s_*$, and absorbing the first term on the right-hand side
into the left-hand side again by choosing $s>0$ large if necessary, we have
\begin{align*}
& se^{2s\mu_0}(\Vert \Delta f\Vert^2_{L^2(\OOO)}
+ \Vert \nabla f\Vert^2_{\{L^2(\OOO)\}^n} + \Vert f\Vert^2_{L^2(\OOO)}) \\
\le & s\int_{\OOO} (\vert \Delta f\vert^2 + \vert \nabla f\vert^2
+ \vert f\vert^2)e^{2s\vadd(x,0)} dx\\
\le & Ce^{2s(\mu_0-r_0)}(\Vert \Delta f\Vert^2_{L^2(\OOO)}
+ \Vert \nabla f\Vert^2_{\{L^2(\OOO)\}^n}) + Ce^{Cs}D^2
\end{align*}
for all $s \ge s_*$.
Noting $C^{-1}\Vert f\Vert^2_{H^2(\OOO)} \le
\Vert \Delta f\Vert^2_{L^2(\OOO)} + \Vert \nabla f\Vert^2_{\{L^2(\OOO)\}^n}
+ \Vert f\Vert^2_{L^2(\OOO)} \le C\Vert f\Vert^2_{H^2(\OOO)}$, we see that
$$
\Vert f\Vert_{H^2(\OOO)} \le Cs^{-1}e^{-2sr_0} \Vert f\Vert^2_{H^2(\OOO)}
+ Ce^{Cs}D^2
$$
for all large $s>s_*$.  Since $r_0 > 0$, for large $s>0$, we can absorb
the first term on the right-hand side into the left -hand side, so that
$$
\Vert f\Vert_{H^2(\OOO)}^2 \le 2Ce^{Cs}D^2.
$$
Thus the proof of Theorem 3 is completed.
\\
\vspace{0.3cm}
\\
{\bf Appendix. Proof of Lemma 4}\\
In order to clarify the dependence in the Carleman estimate
on $\delta_0$ and $M$, for arbitrary $d \in \mathcal{D}$, we prove the
lemma directly and the proof is similar for example to Section 3 in
Yamamoto [21], where a Carleman estimate is proved for a parabolic
equation.

By the usual density argument, it suffices to prove for $f \in C^2_0
(\OOO)$.  We set
$$
\psid(x) = \vad(x,0), \quad
g(x) = e^{s\psid(x)}f(x), \quad q = \Delta f
$$
and
$$
L_0g(x) = e^{s\psid(x)}\Delta(e^{-s\psid(x)}g(x)).
$$
We have
$$
L_0g = \Delta g - 2s\la\psid\nabla d\cdot \nabla g
+ s^2\la^2\psid^2\vert \nabla d\vert^2g - s\la^2\psid\vert\nabla d\vert^2g
- s\la\psid(\Delta d)g \quad \mbox{in $\OOO$}.    \eqno{(1)}
$$
We set
$$
A_1:= - s\la^2\psid\vert\nabla d\vert^2
- s\la\psid(\Delta d).                   \eqno{(2)}
$$
We note that
$$
\vert A_1(x,s,\la)\vert \le Cs\la^2\psi_d \quad \mbox{for $x \in \OOO$ and all
large $\la>0$ and $s>0$}.
$$
Here and henceforth by $C, C_1$, etc., we denote generic constants which are
independent of $s, \la, \psid$ but dependent on $M, \delta_0 > 0$, and
these constants may change line by line.
Moreover we always assume that $\la \ge 1$.  Hence we use the following
inequality: $\la^m \le \la^{m'}$ if $0<m\le m'$.

In view of $A_1$, we have
$$
L_0g = \Delta g - 2s\la\psid\nabla d\cdot \nabla g
+ s^2\la^2\psid^2\vert \nabla d\vert^2g + A_1g
= qe^{s\psid} \quad \mbox{in $\OOO$}.               \eqno{(3)}
$$
Taking into consideration the orders of $(s,\la,\psid)$, we split
$L_0g$ as follows:
$$
\left\{ \begin{array}{rl}
&L_1g = \Delta g + s^2\la^2\psid^2\vert \nabla d\vert^2g + A_1g,\\
&L_2g = - 2s\la\psid\nabla d\cdot \nabla g \quad \mbox{in $\OOO$}.
\end{array}\right.
\eqno{(4)}
$$
Since $\Vert qe^{s\psid}\Vert^2_{L^2(\OOO)}
= \Vert L_1g + L_2g\Vert^2_{L^2(\OOO)}$, we obtain
$$
2\int_{\OOO} (L_1g)(L_2g) dx \le \int_{\OOO} q^2e^{2s\psid} dx.
                                                         \eqno{(5)}
$$
Now we estimate
\begin{align*}
& \int_{\OOO} (L_1g)(L_2g) dx\\
=& -\int_{\OOO} 2s\la\psid (\nabla d\cdot \nabla g) \Delta g dx
- \int_{\OOO} 2s^3\la^3\psid^3\vert \nabla d\vert^2
(\nabla d\cdot \nabla g) g dx\\
\end{align*}
$$
- \int_{\OOO} 2s\la\psid A_1(\nabla d\cdot \nabla g) g dx
=: \sum_{j=1}^3 J_k.           \eqno{(6)}
$$
By the integration by parts and $g=\ppp_{\nu}g = 0$ on $\ppp\OOO$, we have
\begin{align*}
& J_1 = \int_{\OOO} 2s\la\psid \sum_{j,k=1}^n (\ppp_kd)(\ppp_j\ppp_kg)
\ppp_jg dx
+ \int_{\OOO} 2s\la\sum_{j,k=1}^n \ppp_j((\ppp_kd)\psid)(\ppp_jg)\ppp_kg dx\\
=& \int_{\OOO} s\la\psid \sum_{j,k=1}^n (\ppp_kd)\ppp_k(\vert\ppp_jg\vert^2) dx
+ \int_{\OOO} 2s\la \sum_{j,k=1}^n \ppp_j((\ppp_kd)\psid)(\ppp_jg)\ppp_kg dx\\
= & - \int_{\OOO} s\la\sum_{k=1}^n \ppp_k((\ppp_kd)\psid)
\vert \nabla g\vert^2 dx
+ \int_{\OOO} 2s\la\sum_{j,k=1}^n \ppp_j((\ppp_kd)\psid)
(\ppp_jg)\ppp_kg  dx.
\end{align*}
Since
$$
\ppp_j((\ppp_kd)\psid) = (\ppp_j\ppp_kd)\psid + (\ppp_kd)\ppp_j\psid
= \la(\ppp_jd)(\ppp_kd)\psid + (\ppp_j\ppp_kd)\psid,
$$
we have
\begin{align*}
& J_1 = -\int_{\OOO} s\la^2\psid \vert \nabla d\vert^2
\vert \nabla g\vert^2 dx
- \int_{\OOO} s\la(\Delta d)\psid \vert \nabla g\vert^2 dx\\
+& \int_{\OOO} 2s\la^2\psid \sum_{j,k=1}^n (\ppp_jd)(\ppp_kd)(\ppp_jg)
\ppp_kg dx
+ \int_{\OOO} 2s\la \psid\sum_{j,k=1}^n (\ppp_j\ppp_kd)(\ppp_jg)\ppp_kg dx\\
= & -\int_{\OOO} s\la^2\psid \vert \nabla d\vert^2
\vert \nabla g\vert^2 dx
- \int_{\OOO} s\la(\Delta d)\psid \vert \nabla g\vert^2 dx\\
+& \int_{\OOO} 2s\la^2\psid \left\vert \sum_{j,k=1}^n (\ppp_jd)\ppp_jg
\right\vert^2 dx
+ \int_{\OOO} 2s\la \psid\sum_{j,k=1}^n (\ppp_j\ppp_kd)(\ppp_jg)\ppp_kg dx\\
\ge & -\int_{\OOO} s\la^2\psid \vert \nabla d\vert^2
\vert \nabla g\vert^2 dx
- \int_{\OOO} s\la(\Delta d)\psid \vert \nabla g\vert^2 dx\\
+& \int_{\OOO} 2s\la \psid\sum_{j,k=1}^n (\ppp_j\ppp_kd)(\ppp_jg)\ppp_kg dx.
\end{align*}
Therefore
$$
J_1 \ge -\int_{\OOO} s\la^2\psid \vert \nabla d\vert^2
\vert \nabla g\vert^2 dx
- C\int_{\OOO} s\la\psid \vert \nabla g\vert^2 dx.  \eqno{(7)}
$$
Moreover we have
\begin{align*}
& J_2 = -\int_{\OOO} s^3\la^3 \psid^3\vert \nabla d\vert^2
\sum_{k=1}^n (\ppp_kd)\ppp_k(\vert g\vert^2) dx\\
=& \int_{\OOO} s^3\la^3 \sum_{k=1}^n \ppp_k(\psid^3)
\vert \nabla d\vert^2 (\ppp_kd)\vert g\vert^2 dx
+ \int_{\OOO} s^3\la^3\psid^3 \sum_{k=1}^n
\ppp_k(\vert \nabla d\vert^2\ppp_kd)\vert g\vert^2 dx
\end{align*}
$$
= \int_{\OOO} 3s^3\la^4\psid^3 \vert \nabla d\vert^4\vert g\vert^2 dx
- C\int_{\OOO} s^3\la^3\psid^3 \vert g\vert^2 dx       \eqno{(8)}
$$
and
\begin{align*}
& J_3 = -\int_{\OOO} s\la\psid A_1 \sum_{k=1}^n (\ppp_kd)\ppp_k(\vert g\vert^2)
dx\\
=& \int_{\OOO} \sum_{k=1}^n s\la\ppp_k(\psid A_1) (\ppp_kd)\vert g\vert^2 dx
+ \int_{\OOO} \sum_{k=1}^n s\la\psid A_1 (\ppp_k^2d)\vert g\vert^2 dx
\end{align*}
$$
\ge -C\int_{\OOO} s^2\la^4\psid^2 \vert g\vert^2 dx.    \eqno{(9)}
$$
Therefore (6) - (9) yield
\begin{align*}
& \frac{1}{2}\int_{\OOO} \vert q\vert^2 e^{2s\psid} dx
\ge \int_{\OOO} (L_1g)(L_2g) dx \\
\ge & -\int_{\OOO} s\la^2\psid \vert \nabla d\vert^2\vert \nabla g\vert^2 dx
+ \int_{\OOO} 3s^3\la^4\psid^3\vert \nabla d\vert^4 \vert g\vert^2 dx
\end{align*}
$$
- C\int_{\OOO} s\la\psid \vert \nabla g\vert^2 dx
- C\int_{\OOO} (s^3\la^3\psid^3 + s^2\la^4\psid^2)\vert g\vert^2 dx.
                                                                  \eqno{(10)}
$$
For the proof of the lemma, we have to estimate
$s\la^2\psid\vert \nabla g\vert^2 + s^3\la^4\psid^3\vert g\vert^2$,
but the first and the second terms on the right-hand side have
different signs.  Thus we need another estimate.  That is,
multiplying (3) with $-s\la^2\psid\vert \nabla d\vert^2g$ and
integrating over $\OOO$, we obtain
\begin{align*}
& \sum_{j=1}^4 I_k :=
-\int_{\OOO} (\Delta g)s\la^2\psid\vert \nabla d\vert^2 gdx
+ \int_{\OOO} 2s^2\la^3\psid^2 \vert \nabla d\vert^2 (\nabla d\cdot \nabla g)
g dx\\
-& \int_{\OOO} s^3\la^4\psid^3\vert \nabla d\vert^4\vert g\vert^2 dx
- \int_{\OOO} s\la^2\psid A_1\vert \nabla d\vert^2\vert g\vert^2 dx
\end{align*}
$$
= -\int_{\OOO} qe^{s\psid}s\la^2\psid \vert \nabla d\vert^2 g dx.
                                                   \eqno{(11)}
$$
Hence, by integration by parts, we obtain
$$
I_1 = \int_{\OOO} s\la^2\psid \vert \nabla d\vert^2\vert \nabla g\vert^2 dx
+ \int_{\OOO} s\la^2\nabla(\vert \nabla d\vert^2\psid)\cdot (\nabla g)g dx.
$$
Here we have
\begin{align*}
& \int_{\OOO} s\la^2\nabla(\vert \nabla d\vert^2\psid)\cdot (g\nabla g) dx
\le C\int_{\OOO} s\la^3\psid \vert g\vert\vert \nabla g\vert dx\\
=& C\int_{\OOO} (s\la^2\psid \vert g\vert)(\la\vert \nabla g\vert) dx
\le \frac{C}{2}\int_{\OOO} (s^2\la^4\psid^2\vert g\vert^2
+ \la^2\vert \nabla g\vert^2) dx.
\end{align*}
Consequently
$$
I_1 \ge \int_{\OOO} s\la^2\psid \vert \nabla d\vert^2\vert \nabla g\vert^2 dx
- C\int_{\OOO} s^2\la^4\psid^2 \vert g\vert^2 dx
- C\int_{\OOO} \la^2 \vert \nabla g\vert^2 dx.     \eqno{(12)}
$$
Next
\begin{align*}
& I_2 = \int_{\OOO} s^2\la^3\psid^2\vert \nabla d\vert^2
(\nabla d\cdot \nabla(\vert g\vert^2)) dx\\
=& - \int_{\OOO} s^2\la^3 \sum_{k=1}^n \ppp_k(\psid^2)
\vert \nabla d\vert^2(\ppp_kd)\vert g\vert^2 dx
- \int_{\OOO} s^2\la^3\psid^2 \sum_{k=1}^n
\ppp_k((\ppp_kd)\vert \nabla d\vert^2) \vert g\vert^2 dx
\end{align*}
$$
\ge  -C\int_{\OOO} (s^2\la^4\psid^2 + s^2\la^3\psid^2) \vert g\vert^2 dx
\ge  -C\int_{\OOO} s^2\la^4\psid^2 \vert g\vert^2 dx.     \eqno{(13)}
$$
By (2), we see
$$
I_4 \ge -C\int_{\OOO} s^2\la^4\psid^2\vert g\vert^2 dx.   \eqno{(14)}
$$
Since
$$
\left\vert \int_{\OOO} qe^{s\psid} s\la^2\psid \vert \nabla d\vert^2 g dx
\right\vert
\le \frac{1}{2}\int_{\OOO} \vert q\vert^2e^{2s\psid} dx
+ \frac{1}{2}\int_{\OOO} s^2\la^4\psid^2\vert \nabla d\vert^4 \vert g\vert^2
dx,
$$
it follows from (11) - (14) that
\begin{align*}
& \int_{\OOO} s\la^2\psid \vert \nabla d\vert^2 \vert \nabla g\vert^2 dx
- \int_{\OOO} s^3\la^4\psid^3 \vert \nabla d\vert^4 \vert g\vert^2 dx
- C\int_{\OOO} s^2\la^4\psid^2 \vert g\vert^2 dx
- C\int_{\OOO} \la^2\vert \nabla g\vert^2 dx\\
\le& \frac{1}{2}\int_{\OOO} \vert q\vert^2e^{2s\psid} dx
+ \frac{1}{2}\int_{\OOO} s^2\la^4\psid^2 \vert \nabla d\vert^4
\vert g\vert^2 dx.
\end{align*}
Therefore
$$
\int_{\OOO} s\la^2\psid \vert \nabla d\vert^2 \vert \nabla g\vert^2 dx
- \int_{\OOO} s^3\la^4\psid^3 \vert \nabla d\vert^4 \vert g\vert^2 dx
$$
$$
\le \frac{1}{2}\int_{\OOO} \vert q\vert^2 e^{2s\psid} dx
+ C\int_{\OOO} s^2\la^4\psid^2 \vert g\vert^2 dx
+ C\int_{\OOO} \la^2 \vert \nabla g\vert^2 dx.              \eqno{(15)}
$$
Thus we consider (10) $+$ $2\times$ (15):
$$
\int_{\OOO} s\la^2\psid \vert \nabla d\vert^2 \vert \nabla g\vert^2 dx
+ \int_{\OOO} s^3\la^4\psid^3 \vert \nabla d\vert^4 \vert g\vert^2 dx
$$
$$
\le \frac{3}{2}\int_{\OOO} \vert q\vert^2 e^{2s\psid} dx
+ C\int_{\OOO} (s\la\psid + \la^2) \vert \nabla g\vert^2 dx
+ C\int_{\OOO} (s^3\la^3\psid^3 + s^2\la^4\psid^2) \vert g\vert^2 dx.                                                       \eqno{(16)}
$$
Since $\vert \psid(x)\vert \ge e^{-\la M}$ for all $x \in \ooo\OOO$ and
$d\in \mathcal{D}$, we have
$$
C\la^2 = s\la^2\psid \times \frac{C\psid^{-1}}{s}
\le s\la^2\psid \frac{Ce^{\la M}}{s},          \eqno{(17)}
$$
$$
Cs\la\psid = s\la^2\psid \frac{C}{\la}, \quad
Cs^3\la^3\psid^3 = s^3\la^4\psid^3 \frac{C}{\la}     \eqno{(18)}
$$
and
$$
Cs^2\la^4\psid^2 = s^3\la^4\psid^3 \frac{C\psid^{-1}}{s}
\le s^3\la^4\psid^3\frac{Ce^{\la M}}{s}.       \eqno{(19)}
$$
We choose sufficiently large $\la_0 = \la_0(\delta_0,M) > 0$ such that
$\frac{C}{\la_0} \le \frac{1}{4}$ in (18).
For any given $\la \ge \la_0$, we choose $s_1 =s_1(\la,\delta_0,M)$
such that $\frac{Ce^{\la M}}{s_1} \le \frac{1}{4}$ in (17) and (19).
Then
\begin{align*}
&C\la^2, \quad Cs\la\psid \le \frac{1}{4}s\la^2\psid, \\
&Cs^3\la^3\psid^3, \quad Cs^2\la^4\psid^2 \le \frac{1}{4}s^3\la^4\psid^3
\end{align*}
for all $\la \ge \la_0$ and $s \ge s_1$, and noting
$\vert \nabla d\vert \ge \delta_0 > 0$ for all $d \in \mathcal{D}$,
we can absorb the second and the third terms on the right-hand side
of (16) into the left-hand side to obtain
$$
\frac{1}{2}\int_{\OOO} s\la^2\psid \vert \nabla g\vert^2 dx
+ \frac{1}{2} \int_{\OOO} s^3\la^4\psid^3 \vert g\vert^2 dx
\le C\int_{\OOO} \vert q\vert^2 e^{2s\psid} dx.       \eqno{(20)}
$$
We re-write by means of $f$.  Replacing $g = e^{s\psid}f$, we have
$\vert g\vert^2 = \vert f\vert^2 e^{2s\psid}$ and
\begin{align*}
& \vert \nabla f\vert^2 e^{2s\psid}
= \vert \nabla g - s\la(\nabla d)\psid e^{s\psid}f\vert^2\\
\le & 2\vert \nabla g \vert^2
+ 2s^2\la^2\vert \nabla d\vert^2\psid^2 \vert f\vert^2 e^{2s\psid}
\le 2\vert \nabla g \vert^2
+ Cs^2\la^2\psid^2 \vert f\vert^2 e^{2s\psid}.
\end{align*}
Substituting them into (20), we complete the proof of
Lemma 4.


\begin{thebibliography}{99} %
\bibitem{A} A. Kh. Amirov, Integral Geometry and Inverse Problems for
Kinetic Equations, VSP, Utrect, 2001.

\bibitem{Ba1} G. Bal,
Inverse transport theory and applications, Inverse Problems {\bf 25} (2009),
053001, 48 pp.

\bibitem{Bo}
L. Baudouin, M. De Buhan and S. Ervedoza,
Global Carleman estimates for waves and applications,
Comm. Par. Differ. Equ. {\bf 38} (2013), 823-859.

\bibitem{Bei} L. Beilina and M.V. Klibanov,
Approximate Global Convergence and Adaptivity for Coefficient
Inverse Problems, Springer-Verlag, Berlin, 2012.

\bibitem{Bel} S.P. Belinskij, On one of the inverse problems for
linear symmetrical $t-$hyperbolic systems with $n+1$ independent variables,
Diff. Equ. {\bf 12} (1976), 15-23.

\bibitem{BM}
M. Bellassoued and M. Yamamoto,
Carleman Estimates and Applications to Inverse Problems for
Hyperbolic Systems, Springer-Verlag, Berlin, 2015.

\bibitem{BK}
A.L. Bugheim and M.V. Klibanov,
Global uniqueness of class of multidimensional inverse problems,
Soviet Math. Dokl. {\bf 24} (1981), 244-247.

\bibitem{Case}
K.M. Case and P.F. Zweifel, Linear Transport Theory,
Addison-Wesley, Reading (MA), 1967.

\bibitem{GO} P. Gaitan and H. Ouzzane,
Inverse problem for a free transport equation using Carleman
estimates, Applicable Analysis {\bf 93} (2014), 1073-1086.

\bibitem{Ho}
L. H\"ormander, Linear Partial Differential Operators, Springer-Verlag,
Berlin, 1963.

\bibitem{IY1}
O.Yu.Imanuvilov and M.Yamamoto, Global Lipschitz stability in an
inverse hyperbolic problem by interior observations,
Inverse Problems {\bf 17} (2001), 717-728.

\bibitem{IY2}
O.Yu. Imanuvilov and M. Yamamoto, Global uniqueness and stability in determining
coefficients of wave equations,
Comm. Par. Diff. Equ.{\bf 26} (2001), 1409-1425.

\bibitem{Is} V. Isakov,
Inverse Problems for Partial Differential Equations, Springer-Verlag,
Berlin, 2006.

\bibitem{Kli}
M.V. Klibanov, Inverse problems and Carleman estimates,
Inverse Problems {\bf 8} (1992), 575-596.

\bibitem{KliPa}
M.V. Klibanov and S.E. Pamyatnykh,
Global uniqueness for a coefficient inverse problem for the
non-stationary transport equation via Carleman estimate,
J Math. Anal. Appl. {\bf 343} (2008), 352-365.

\bibitem{KY}
M.V. Klibanov and M. Yamamoto,
Exact controllability for the non stationary transport equation,
SIAM Control Optim. {\bf 46} (2007), 2071-2195.

\bibitem{MY}
M. Machida and M. Yamamoto,
Global Lipschitz stability in determining coefficients
of the radiative transport equation, Inverse Problems {\bf 30}
(2014), 035010.

\bibitem{Ba2} K. Ren, G. Bal and A.H. Hielscher,
Transport- and diffusion-based optical tomography in small
domains: a comparative study, Applied Optics {\bf 27} (2007), 6669-6679.


\bibitem{R}
V.G. Romanov, Inverse Problems of Mathematical Physics,
VNU Science Press, Utrech, 1987.

\bibitem{S}
P. Stefanov, Inverse problems in transport theory, in "Inside out:
inverse problems and applications", Cambridge University Press,
Cambridge, pp.111-131, 2003.

\bibitem{Y}
M. Yamamoto, Carleman estimates for parabolic equations and applications,
Inverse Problems {\bf 25} (2009), 123013 (75pp).

\end{thebibliography}
\end{document}